\documentclass[12pt]{article}

\usepackage{amsmath,amsthm}
\usepackage{amssymb}
\usepackage{amsfonts}
\usepackage{amscd}

\usepackage{color}

\usepackage[normalem]{ulem}

\usepackage{cancel}

\newtheorem{thm}{Theorem}[section]
\newtheorem{theorem}[thm]{Theorem}

\newtheorem{lemma}[thm]{Lemma}
\newtheorem{proposition}[thm]{Proposition}

\theoremstyle{definition}
\newtheorem{definition}[thm]{Definition}

\newtheorem{remark}[thm]{Remark}
\newtheorem{observation}[thm]{Observation}
\newtheorem{free text}[thm]{}

\newcommand{\bG} {G} 
\newcommand{\Z} {\mathbb{Z}}  
\newcommand{\C} {\mathbb{C}}  

\newcommand{\bk} {k}
\newcommand{\fg} {\mathfrak {g}}
\newcommand{\fgl} {\mathfrak {gl}}
\newcommand{\fh} {\mathfrak {h}}

\newcommand{\salg} {\hbox{\rm (salg)}}

\newcommand{\lie} {\text{(Lie)}}

\newcommand{\sets}{{(\hbox{sets})}}

\newcommand{\chesgrp} {\hbox{\rm (chesgrps)}}
\newcommand{\CSHCP} {\hbox{\rm (CSHCP)}}

\newcommand{\lra} {\longrightarrow}

\newcommand{\rgl}{\mathfrak{gl}}

\newcommand{\rGL}{\mathrm{GL}}

\newcommand{\Lie}{\hbox{\sl Lie}}

\newcommand{\al}{\alpha}
\newcommand{\be}{\beta}
\newcommand{\ga}{\gamma}

\newcommand{\cO}{\mathcal{O}}

\newcommand{\U} {\mathcal{U}}
\newcommand{\kzg} {\mathcal{U}_\Z(\fg)}
\newcommand{\cV}{{\cal {V}}}
\begin{document}

\bigskip

\centerline{\LARGE \bf
Algebraic supergroups with}

\medskip

\centerline{\LARGE \bf Lie superalgebras of classical type}

\vskip45pt

\centerline{Rita Fioresi$^\flat$, Fabio Gavarini$^\#$}

\bigskip

\centerline{\it $^\flat$ Dipartimento di Matematica, Universit\`a
di Bologna }
\centerline{\it Piazza di Porta San Donato, 5  ---
I--40127 Bologna, ITALY}
\centerline{{\footnotesize e-mail: fioresi@dm.unibo.it}}

\bigskip

\centerline{\it $^\#$ Dipartimento di Matematica,
Universit\`a di Roma ``Tor Vergata'' }
\centerline{\it via della ricerca scientifica 1  ---
I--00133 Roma, ITALY}
\centerline{{\footnotesize e-mail: gavarini@mat.uniroma2.it}}

\vskip65pt

\begin{abstract}
 We show that every connected affine algebraic supergroup defined over a field $k$, with diagonalizable maximal torus and whose tangent Lie superalgebra is a  $ k $--form  of a complex simple Lie superalgebra
of classical type is a Chevalley supergroup, as it is defined and constructed explicitly in  [R.~Fioresi, F. Gavarini,  {\it Chevalley Supergroups},  Memoirs of the AMS  {\bf 215}
 (2012), no.~1014].\footnote{\ 2000 {\it MSC}\;: \, Primary 14M30, 14A22; Secondary 58A50, 17B50.}
\end{abstract}

\vskip45pt

\section{Introduction}
 In \cite{fg1} we have given the supergeometric analogue of the classical Chevalley's construction (see \cite{st}), which enabled us to build a supergroup out of data involving only a complex Lie superalgebra  $ \fg $ of classical type and a suitable complex faithful representation.  Such a supergroup is affine connected, with associated classical subgroup being reductive  $ k $--split  (i.e.~it admits a diagonalizable maximal torus) and with tangent Lie superalgebra isomorphic to  $ \fg \, $:  thus we obtained an  {\sl existence result\/}  for such supergroups.  In particular, this provided the first unified construction of affine algebraic supergroups with tangent Lie superalgebras of classical type; in particular, it was also (as far as we know) the very first explicit construction of algebraic supergroups corresponding to the simple Lie superalgebras of basic exceptional type.

\smallskip

   In this paper we tackle the  {\sl uniqueness\/}  problem, cast in the following form: ``is any such supergroup isomorphic to a supergroup obtained via the Chevalley's construction''?  Our answer is positive.

\smallskip

   We start with an affine algebraic supergroup  $ G \, $,  defined over a field  $ k \, $ with associated classical subgroup  $ G_0 $  which is  $ k $--split  reductive, and with tangent Lie superalgebra a  $ k $--form  of a complex Lie superalgebra of classical type (plus a consistency condition): then we prove that  $ G $  is given by our Chevalley supergroup construction.  Note that all the conditions we impose actually are necessary, as they do hold for Chevalley supergroups.

\smallskip

   As  $ G_0 $  is  $ k $--split  reductive, by Chevalley-Demazure theory it can be realized via the Chevalley construction as a closed subgroup of some  $\rGL\big(\widetilde{V}\big) \, $,  where  $ \widetilde{V} $  is a suitable  $ G_0 $--module.  Let  $ \widetilde{V}^* $  be the dual  $ G_0 $--module.  Since  $ G $  is an affine supergroup over a field  $ k \, $,  it is linearizable, that is  $ \, G \subseteq \rGL_{m|n} \, $  (for suitable  $ m $  and  $ n $),  hence we can build the induced  $ {(\rGL_{m|n})}_0 $--module  $ \, U := \text{\it Ind}_{\,\bG_0}^{\;{(\rGL_{m|n})}_0} \big(\widetilde{V}^*\big) \, $  and its dual  $ \, U^* $,  which both are naturally  $ {(\rgl_{m|n})}_0 $--modules  as well: note also that  $ U^* $ contains a  $ G_0 $--submodule  isomorphic to  $ \widetilde{V} \, $.  Inducing then for the Lie superalgebras we get the $ \rgl_{m|n} $--module $ \, W := \text{\it Ind}_{\,{(\rgl_{m|n})}_0}^{\;\;\rgl_{m|n}}\!\big(U^*\big) \, = \, \U\big(\rgl_{m|n}\big) \otimes_{{}_{\scriptstyle \U((\rgl_{m|n})_0)}} \! U^* \, $.  Now  $ W $ is also a $ \rGL_{m|n} $--module  and (by restriction) a  $ G $--module:  moreover, it contains the (finite-dimensional)  $ G $--submodule  $ \, V := \U(\fg) \otimes_{\U(\fg_0)} \! \widetilde{V} \, $,  where  $ \widetilde{V} $  is identified with a  $ G_0 $--submodule of  $ U^* $. {\sl N.B.:}  for the sake of simplicity of exposition, we are hiding here several technicalities, to be specified later on in the main text.

\smallskip

   The very construction of  $ V $  allows us to build the Chevalley supergroup  $ G_V $  associated with the  $ \fg $--representation  $ V $  and to view both  $ G $  and  $ G_V $  as closed subgroups of the  {\it same\/}  $ \rGL(V) \, $.  The last step is to note that both  $ G $  and  $ G_V $  are globally split   --- as any affine supergroup over a field, by  Theorem 4.5 in  \cite{ms1}.  Since the ordinary algebraic groups are the same,  $ \, G_0 = {(G_V)}_0 \, $,  we have that both supergroups are smooth as well.  We conclude then  $ \, G = G_V \, $  by infinitesimal considerations, since they have the same Lie superalgebra.

\smallskip

   In the last section we make some important remarks between the
equivalence of categories of certain Super Harish-Chandra pairs
and the algebraic supergroups we have studied in the present work.

\smallskip

   Parallel constructions and results, concerning existence (by a Chevalley like construction) and uniqueness of algebraic supergroups associated with simple Lie superalgebras of Cartan type are presented in \cite{ga}.

\vskip17pt

   \centerline{\bf Acknowledgements}
 \vskip3pt
   We thank prof.~V.~S.~Varadarajan, whose valuable suggestions on
a previous version of the manuscript helped us to improve our work.
We also thank the referee for his/her careful and deep analysis of
the present paper.

\bigskip

\section{Chevalley supergroups} \label{chevalley}

In this section we review briefly the construction of Chevalley supergroups (see \cite{fg1}, \cite{fg2}) and then we discuss some of their properties.  For all details about the construction we refer to \cite{fg1}.  The new property that we present here is that every Chevalley supergroup  $ G_V \, $,  defined as a subgroup of some  $ \rGL(V) \, $,  is in fact  {\sl closed\/}  in  $ \rGL(V) \, $.

\medskip

\subsection{Definition of Chevalley supergroups}  \label{def_che-sgroups}

Let $\fg$ be a complex Lie superalgebra of classical type and $\fh$
a fixed Cartan subalgebra of  $ \fg_0 \, $.  Then we have the corresponding root system  $ \, \Delta = \Delta_0 \cup \Delta_1 \, $,  with  $ \Delta_0 $ and  $ \Delta_1 $  being the sets of even and
of odd roots respectively: these roots are the non-zero eigenvalues of the (adjoint) action of  $ \fh $ on  $ \fg \, $,  while the corresponding eigenspaces, resp.~eigenvectors, are called  {\it root spaces},  resp.~{\it root vectors}.  For root vectors, we adopt the simplified notation of the cases when  $ \fg $  is not of type  $ A(1,1) $,  $ P(3) $  or  $ Q(n) $   --- cf.~\cite{ka} ---   but  {\sl all what follows holds for those cases too},
%
%
and all our results hold for all complex Lie superalgebras of classical type, but
for the cases  $ D(2,1;a) $  when  $ \, a \not\in \Z \, $.
                                                    \par
   Like in the classical setting, one can define special elements  $ \, H_\alpha \in \fh \, $,  called  {\it coroots},  associated with the roots  $ \alpha \, $.

\medskip

   A key notion in  \cite{fg1}  is that of  {\it Chevalley basis\/}  of  $ \fg \, $.  This is any  $ \C $--basis  of $ \fg $  of the form
  $$  B  \; = \;  \big\{ H_1 \dots H_\ell \big\} \cup \big\{\, X_\alpha \, , \, \alpha \in \Delta \,\big\} $$
such that  (cf.~\cite{fg1},  Def.~3.3):
\begin{itemize}
\item the  $ H_i $'s, called  the {\sl Cartan elements\/}  of  $ B \, $,
form a  $ \C $--basis  of  $ \fh $  (with some additional properties);
\item every  $ X_\alpha $  is a root vector associated with the root
$ \alpha \, $;
\item  the structure coefficients for the Lie superbracket in  $ \fg $
with respect to these basis elements are integers with some special
properties.
\end{itemize}

\smallskip

   The very existence of Chevalley bases is proved in \cite{fg1}, sec.~3.

\medskip

   If  $ B $  is a Chevalley basis of  $ \fg $  as above, we set  $ \;
\fg_\Z := \text{\it Span}_\Z\{B\} \; \big(\! \subseteq \fg \big) \; $  for its  $ \Z $--span.  Moreover, we define an important integral lattice inside  $ \U(\fg) \, $,  namely the {\it Kostant superalgebra}.  This is the  $ \Z $--supersubalgebra  $ \kzg $  of  $ \U(\fg) \, $  generated by the following elements: all divided powers in the even root vectors of  $ B \, $,  all odd root vectors of  $ B \, $,  and all binomial coefficients in the Cartan elements of  $ B \, $  (see  \cite{fg1}, sec.~4.1).

\medskip

   We associate to  $ \kzg $  the notion of  {\sl admissible lattice\/}  in a  $ \fg $--module:

\medskip

\begin{definition} \label{rational}
Let  $ \fg \, $,  $ \, B = \big\{ H_1 \dots H_\ell \big\} \cup \big\{\, X_\alpha \, , \alpha \in \Delta \,\big\} \, $  and  $ \kzg $  be as above.  Let  $ V $  be a complex finite dimensional  $ \fg $--module.
 We say that $V$ is  {\it rational\/}  if the  $ H_i $'s  act diagonally on  $ V $  with integral eigenvalues.
 We say that an integral lattice  $ M $ in $V$   --- that is, a free  $ \Z $--submodule  $ M $  of  $ V $ such that $ \, \textit{rk}_{\,\Z}(M) = \textit{dim}_{\,\C}(V) $  ---   is
{\it admissible\/}  if it is  $ \kzg $--stable.
\end{definition}

\medskip

   Given a complex representation $V$ of  $ \fg $  as above, there exists always an admissible lattice  $ M $  and an integral form  $ \fg_V $  of  $ \fg $  keeping such a lattice stable (see  \cite{fg1}, \S 5.1).  This allows us to shift from the complex field  $ \C $  to any commutative unital ring  $ k \, $.

\medskip

\begin{definition} \label{crucialdef}
 Let the notation be as above, and assume also that the representation  $ V $  is  {\sl faithful}.  For any fixed commutative unital ring  $ k \, $,  define
 \vskip4pt
   \centerline{ $ \fg_k \, := \, k \otimes_\Z \fg_V \; ,
\qquad  V_k := k \otimes_\Z M \; ,
\qquad \U_k(\fg) := k \otimes_\Z \kzg $ }
 \vskip4pt
 Then we say that  $ \fg_k \, $,  resp.~$ M $,  is a \textit{$k$--form} of  $ \fg \, $,  resp.~of  $ V_k $.
\end{definition}

\smallskip

\begin{remark}  \label{super-dist}
 For any algebraic supergroup  $ G $,  one can introduce the notion of  {\sl superalgebra of distributions}  $ \text{\sl Dist}_k(G) \, $,  by an obvious extension of the standard notion in the even setting; see  \cite{bk},  \S 4, for details.  One can easily see   --- like in  \cite{bk}, \S 4  ---   that  $ \; \text{\sl Dist}_k(G) = \U_k(\fg) \; $;  in particular, this shows that  $ \U_k(\fg) $  is independent of the choice of a specific Chevalley basis in  $ \fg \, $.
                                                \par
   More important (for later use), is the fact that if  $ \, \varphi : G' \longrightarrow G'' \, $  is a morphism between two supergroups, then it induces (functorially) a morphism  $ \, D_\varphi : \text{\sl Dist}_k\big(G'\big) \longrightarrow \text{\sl Dist}_k\big(G''\big) \, $,  which is injective whenever  $ \varphi $  is injective.  If in addition  $ G'$  and  $ G'' $  satisfy the assumptions we gave above for  $ G \, $,  so that  $ \, \U_k\big(\fg'\big) = \text{\sl Dist}_k\big(G'\big) \, $  and  $ \, \U_k\big(\fg''\big) = \text{\sl Dist}_k\big(G''\big) \, $,  we have then  $ \, D_\varphi \! : \U_k\big(\fg'\big) \! \rightarrow \U_k\big(\fg''\big) \, $,  which is an embedding if  $ G' $  is subsupersgroup of  $ G'' $.
\end{remark}

\medskip

   We need now to recall the notion of commutative superalgebras.
                                                                \par
   We call  {\it  $ k $--superalgebra\/}  any associative, unital  $ k $--algebra  $ A $  which is  $ \Z_2 $--graded  (as a  $ k $--algebra):  so  $ A $  bears a  $ \Z_2 $--splitting  $ \, A = A_0 \oplus A_1 \, $  into direct sum of super-subvector spaces, with  $ \, A_a \, A_b \subseteq A_{a+b} \, $.  We define the  {\it parity\/}  $ \, |a| \in \Z_2 \, $  of any  $ \, a \in \big( A_0 \cup A_1 \big) \setminus \{0\} \, $  by the condition  $ \, a \in A_{|a|} \, $;  the elements in  $ A_0 $  are called  {\it even},  those in  $ A_1 $  {\it odd}.  All  $ k $--superalgebras  form a category, whose morphisms are those in the category of  $ k $--algebras  which preserve the unit and the  $ \Z_2 $--grading.
                                                                \par
   A  $ k $--superalgebra  $ A $  is said to be  {\it commutative\/}  iff  $ \; x \, y = (-1)^{|x|\,|y|} y \, x \; $  for all homogeneous  $ \, x $,  $ y \in A \, $  and  $ \; z^2 = 0 \; $  for all odd  $ \, z \in A_1 \, $.  We denote  by $ \salg $   --- or $ \salg_k $  ---   the category of commutative  $ k $--superalgebras.

\vskip4pt

   As a matter of notation, we write  $ \text{(grps)} $  for the category of groups.

\medskip

   Finally, we are ready to give the definition of Chevalley supergroup over the commutative ring  $ k \, $.

\medskip

\begin{definition} \label{chevgroup}
 Let the notation be as above.  We define  {\it Chevalley supergroup\/}
the supergroup functor  $ \; \bG_V : \salg_\bk \longrightarrow \text{(grps)} \; $
defined as:
 $ \; \bG_V(A) \, := \,
\Big\langle\, {\bG_{V,0}(A)} \, , \, 1+\theta_\beta X_\beta \,\;
\Big|\;\, \beta \in \Delta_1, \, \theta_\be \in A_1 \Big\rangle \; $  $ \; \Big(\! \subseteq \rGL\big(V_k\big)(A) \Big) \; $,
 for all $ \, A \in \salg_\bk \, $,  where $\bG_{V,0}$ is the ordinary reductive group scheme
%
%
 associated via the Chevalley recipe with the  $ \bG_{V,0} $--module  $ V_k $  (cf.~\cite{fg1}, sec.~5).  As usual $\rGL(V_k)$ denotes the general linear supergroup scheme
\end{definition}

\bigskip

Let us fix a total order
(with some mild conditions) in  $ \Delta_1 \, $,  and let
$\bG_{V,1}^<$  be the functor of points of the superscheme corresponding to
ordered products of elements of the type  $ \, 1+\theta X \in \bG_V(A) \, $
where $X$ is a positive root vector.
We have that  $ \, \bG_{V,1}^< \cong \mathbb{A}^{0|N} \, $  where
$ \, N = \textit{dim}_{\,\C}(\fg_1) = \big| \Delta_1 \big| \, $  and
$\mathbb{A}^{0|N}$ denotes the purely odd affine superspace
(see \cite{fg1}, sec.~5, and \cite{fg2}, sec.~4 for details).

\bigskip

\begin{theorem} \label{mainreschev}
 The group product  $ \, \bG_{V,0} \times \bG_{V,1}^< \lra \bG_V \, $
induces an isomorphism of superschemes.  In particular we have
 $ \; \bG_V \, \cong \, \bG_{V,0} \times \mathbb{A}^{0|N} \; $
(with  $ N $  as above), so that  $ \bG_V $
is an affine supergroup scheme (it is representable).
\end{theorem}

   Theorem  \ref{mainreschev}  is the main result in  \cite{fg1}:  in particular, it states the representability of the supergroup functor  $ \bG_V \, $,  so that the terminology  {\sl Chevalley supergroup\/}  is fully justified.  Furthermore, for  $ k $  a field we have  $ \, \Lie(\bG_V) = \fg_k \, $  as expected.  Finally since by the classical theory  $ \bG_{V,0} $  is  {\sl connected},  $ G_V $  is connected.

%
%
%
%
%

\bigskip

\subsection{The Chevalley supergroup  $ \bG $  is closed inside  $ \rGL(V_k) \, $}

Let $k$ be a unital commutative ring.
All our algebras and modules will now be over $k$ unless
otherwise specified.

\medskip

We now wish to prove that $\bG_V$ embeds naturally into the general linear supergroup
$\rGL(V_k)$ as a {\sl closed subsuperscheme}.  Note that, when  $ k $  is a field, the affine supergroup  $\bG_V$  embeds into some
$\rGL(W)$ as a closed supergroup subscheme (see \cite{ccf}, ch.~11);
we now want to show that we can always choose
$ \, W := V_k \, $,  where  $ V_k $  is the
$ \fg $--supermodule  used to construct  $ \bG_V $  itself.

\medskip

Let us start with some observations.

\medskip

Let $\fgl(V_k)$ be the Lie superalgebra of all the endomorphisms
of the free module  $ V_k \, $:  we denote with $\fgl(V_k)_0$  the set of all the endomorphisms
preserving parity, and with $\fgl(V_k)_1$ the set of those reversing parity.
Its functor of points  $ \, \fgl(V_k) : \salg \lra \lie \, $ is Lie algebra valued
(hereafter $\lie$ denotes the category of Lie algebras) and it is given by:
$$
 \fgl(V_k)(A) \; := \; {\big( A \otimes \fgl(V_k) \big)}_0 \; = \;
A_0 \otimes \fgl(V_k)_0 \, \oplus \, A_1 \otimes \fgl(V_k)_1
$$
Notice that in this equality the symbol $\fgl(V_k)$ appears with
two very different meanings: on the left hand side it is a
Lie algebra valued functor, while on the right hand side it is
just a free module over $k$. This is a most common abuse
of notation in the literature.
Hence  $\fgl(V_k)(A)$  splits into direct sum of
$$
\fgl(V_k)_0(A) \, = \, A_0 \otimes \fgl(V_k)_0 \;\; ,  \qquad
\fgl(V_k)_1(A) \, = \, A_1 \otimes \fgl(V_k)_1
$$
corresponding respectively to the functor of points of the
purely even Lie superalgebra $\fgl(V_k)_0$   --- hence a Lie algebra ---
and to the functor of points of the purely odd superspace
$ \fgl(V_k)_1 \, $.
Now define the functor  $ \, \rGL(V_k)_1 : \salg \lra \sets \, $  by
$$
\rGL(V_k)_1(A) \, = \, I + \fgl(V_k)_1(A)   \eqno \forall \;\; A \in \salg  \qquad
$$
where $I$ denotes the identity in  $ \rGL(V_k)_1(A) \, $.
One can check immediately that this is a representable functor
corresponding to the affine purely odd superspace  $ \mathbb{A}^{0|2mn} $,  where
$m|n$ is the dimension of  $ V_k \, $.
One also sees easily that $\rGL(V_k)_1$ is a subfunctor
and a subscheme of  $ \rGL(V_k) \, $.  The reader must be warned that
$\rGL(V_k)_1$ has no natural supergroup structure.

\medskip

The next proposition clarifies the relation between
$\rGL(V_k)_1$ and  $ \rGL(V_k) \, $.

\medskip

\begin{proposition} \label{GLdec}
 Let the notation be as above.
    Then the multiplication map\break
 $ \; \rGL(V_k)_0 \times \rGL(V_k)_1 \relbar\joinrel\lra \rGL(V_k) \; $  induces an isomorphism of superschemes, where $\rGL(V_k)_0$ denotes as usual
the closed superscheme of $\rGL(V_k)$ corresponding to the ordinary
underlying affine group.  In particular, both   $\rGL(V_k)_0$ and $\rGL(V_k)_1$
are closed supersubschemes of  $ \rGL(V_k) \, $.
\end{proposition}

\begin{proof}
Given  $ \, A \in \salg \, $,  let us consider an  $ A $--point  of  $ \rGL(V_k) \, $,  say
$$
\left(\begin{array}{cc} a & \be \\ \ga & d \end{array}\right) \in \rGL(V_k)(A)
$$
Then  $ \, a \, $, $d$ $\in \rGL(V_k)_0 \, $  are invertible matrices and this
immediately allows us to build the inverse morphism of the map
$ \, \rGL(V_k)_0 \times \rGL(V_k)_1 \lra \rGL(V_k) \, $  given by restriction of the multiplication, namely
$$
\begin{array}{ccc}
\rGL(V_k) & {\buildrel \cong \over
{\lhook\joinrel\relbar\joinrel\relbar\joinrel\relbar\joinrel\relbar\joinrel\twoheadrightarrow}} &
\hskip-25pt
 {\rGL(V_k)}_0 \times {\rGL(V_k)}_1 \\ \\
\Big(\! {{{\;\; a \;|\; \beta \;\,}
\over {\,\; \gamma \;|\; d \;\,}}} \!\Big) & \mapsto &
 \hskip-11pt
\bigg( \Big(\! {{{\,\; a \;|\; 0 \;\,} \over {\,\; 0 \;|\; d \;\,}}} \!\Big) \, , \,
\Big(\! {{{\,\; \hskip6pt I_{m} \hskip3pt \;|\; a^{-1} \beta \;\,}
\over {\,\; d^{-1} \gamma \;|\; \hskip6pt I_{n} \hskip4pt \;\,}}} \!\Big) \bigg)
\;
\end{array}
$$
where $m|n$ is the dimension of $V_k$ and $I_s$ is the identity matrix of size  $ s \, $. \\
The statement about $\rGL(V_k)_0$ and $\rGL(V_k)_1$ being closed is clear.
\end{proof}

\medskip

\begin{theorem} \label{GVclosed}
%
%
 Let  $ \, \bG_V $  be the Chevalley supergroup associated with
the complex Lie superalgebra $\fg$ and to a complex
representation  $ V $ of\/  $ \fg \, $.
   \hbox{Then  $ \, \bG_V $  is}
 a closed supergroup subscheme
in the general linear supergroup scheme  $ \rGL(V_k) $.
\end{theorem}

\begin{proof}
 By the very definition of Chevalley supergroup and by
Theorem \ref{mainreschev}
we have that
$$
\bG_V \; \cong \; \bG_{V,0} \times \bG_{V,1}^< \; \subseteq \;
\rGL(V_k) \; \cong \; \rGL(V_k)_0 \times \rGL(V_k)_1
$$
By the classical theory we have that $\bG_{V,0}$ is a closed subgroup (scheme) of  $\rGL(V_k)_0 \, $,  thus it is enough to show that  $ \, \bG_{V,1}^< \, $  is closed too   --- as a super-subscheme of  $ \rGL(V_k) \, $.
 \vskip3pt
   Let us look closely at the embedding of  $ \, \bG_{V,1}^< \, $  inside  $ \, \rGL(V_k) \, $.
By  Theorem \ref{mainreschev}
we have an isomorphism
$ \, \Psi : \mathbb{A}^{0|N} \! \lra G_{V,1}^< \; $  given by
 \vskip3pt
   \centerline{ $  \Psi_A \, : \, \mathbb{A}^{0|N}(A) \! \lra G_{V,1}^<(A)  \quad ,  \qquad  (\vartheta_1,\dots,\vartheta_N) \, \mapsto \, {\textstyle \prod_{i=1}^N} \,
x_{\gamma_i}(\vartheta_i) $ }
 \vskip3pt
\noindent
 where the product in right-hand side is ordered w.r.t.~some total order on  $ \Delta_1 $  for which  $ \Delta_1^+ $  follows  $ \Delta_1^- \, $,  or viceversa.  In particular, the point  $ 0 $  in  $ \mathbb{A}^{0|N} $  corresponds to the identity  $ I $  in  $ \bG_{V,1}^< \, $;  thus the tangent superspace to  $ \bG_{V,1}^< $  at  $ I $  corresponds to the tangent superspace to  $ \mathbb{A}^{0|N} $  at  $ 0 \, $,  naturally identified with  $ \mathbb{A}^{0|N} $  again.
                                                                     \par
   Given  $ \, A \in \salg \, $,  we have for
$ \; g = {\textstyle \prod_{i=1}^N} \, x_{\gamma_i}(\vartheta_i) \in \bG_{V,1}^<(A) \; $:
$$
g  \, = \,  {\textstyle \prod_{i=1}^N} \, x_{\gamma_i}(\vartheta_i)  \, = \,
I + {\textstyle \sum_{i=1}^N} \, \vartheta_i \, X_{\gamma_i} + \cO(2)  \;
\in \;  \mathfrak{gl}\big(V_k(A)\big)   \eqno (\star)
$$
where  $ \, \cO(2) \, $  is some element in  $ \, \mathfrak{gl}\big(V_k(A)\big)
= A_0 \otimes_k \mathfrak{gl}(V_k)_0 + A_1 \otimes_k \mathfrak{gl}(V_k)_1 \, $
whose (non-zero) coefficients in  $ A_0 $  and  $ A_1 $  actually belong to
$ J_A^{\,2} \, $,  the ideal of  $ A $
generated by  $ \, A_1^{\,2} := A_1 \cdot A_1 \, $.

Consider now the  closed subscheme $H$ in $\rGL(V_k)_1$
whose functor of points is defined as
  $$  H(A)  \; := \;  I \, + \, {\textstyle \sum_i} \, \vartheta_i \, X_{\gamma_i}  $$
We have an invertible natural transformation $\phi$
$$
\begin{array}{ccccl}
\phi_A: & G_{V,1}^<(A) & \lra & H(A) & \Big(\! \subseteq \rGL(V_k)(A) \Big)  \\ \\
& {\textstyle \prod_{i=1}^N} \, x_{\gamma_i}(\vartheta_i) & \mapsto &
I+ \sum_i \vartheta_i X_{\gamma_i} &
\end{array}
$$
which maps $G_{V,1}^<$ isomorphically onto the closed subscheme  $ H $ in  $ \rGL(V_k)_1 \, $,
whence $G_{V,1}^<$ is a closed subsuperscheme of  $ \rGL(V_k)_1 \, $.
\end{proof}

\section{Uniqueness Theorem}  \label{unique}

   Hereafter, we assume $ k $  to be a  {\sl field},  with  $ \, {\mathrm{char}}(k) \neq 2 \, , 3 \, $.

\smallskip

   In this section we prove the main result of our paper, which we summarize
as follows. Let  $ G $  be a connected affine algebraic supergroup, whose
tangent Lie superalgebra  $ \fg_k $  is a  $ k $--form  of a complex Lie
superalgebra of classical type (see  Def.~\ref{crucialdef});  we assume
also that its even subgroup  $ \bG_0 $  is reductive and  $ k $--split,
i.e.~it admits a diagonalizable maximal torus. We assume further that
$ (\fg_k)_0 $, the even part of $\fg_k$ is an ingredient in the
recipe that allows us to realize the ordinary group $ G_0 \, $ as
a Chevalley group.

\medskip

We then show that such a $ G $  is isomorphic to a Chevalley supergroup
$ G_V $  as we constructed in  \cite {fg1} according to the recipe described
in the previous section.

\medskip

  We start with a result relative to the chosen admissible representation $V$ of  the complex Lie superalgebra  $ \fg \, $,  inducing the embedding of  $ G_V $  in  $ \rGL(V_k) \, $.

\subsection{Linearizing $\bG$}  \label{linearizing-bG}

   Let  $ \bG $  be a connected affine algebraic supergroup over
$ k \, $ and let  $ \, \fg_k := \Lie\,(\bG) \, $  be the tangent
Lie superalgebra of  $ G \, $.

\medskip

We assume $\fg_k$ to be a  $ k $--form  of a complex Lie superalgebra
$ \fg \, $,  that is  $ \, \fg_k = k \otimes \fg^\Z \, $
(cf.~Definition \ref{crucialdef}), where here  $ \fg^\Z $
is any integral lattice inside
the complex Lie superalgebra  $ \fg \, $.  Moreover, we assume
the complex Lie superalgebra  $ \fg $  to be simple of
{\sl classical type\/}  (in the sense of Kac's terminology, see  \cite{ka}).
It follows that the even part  $ \fg_0 $  of  $ \fg $  is a
{\sl reductive\/}  Lie algebra.
 \vskip4pt
   Let  $ G_0 $  be the ordinary subgroup underlying  $ G \, $:  its tangent Lie algebra is  $ \, \Lie(G_0) = {\Lie(G)}_0 = {(\fg_k)}_0 \, $.  We assume that $G_0$ is reductive and  $ k $--{\sl split},  i.e.~it admits a diagonalizable maximal torus.

\medskip

By the classical theory then  $ \bG_0 $  can be realized via the classical Chevalley construction (see for example \cite{ja2}, part II, 1.1).  In short, there exists a complex  $ \fg_0 \, $--module  $ \widetilde{V} $  which is faithful, rational, finite-dimensional, so that  $\bG_0 $  is isomorphic to the affine group-scheme (over  $ \Z $)  associated with  $ \fg_0 $  and  $ \widetilde{V} $  by the classical Chevalley's construction (see also Demazure \cite{dem}),  using some admissible lattice  $ \widetilde{M} $  in  $ \widetilde{V} $.  Here such words as  {\sl rational\/}  and  {\sl admissible\/}  refer to the choice of any Chevalley basis  $ B'_0 $  (in the classical sense) of the reductive Lie algebra  $ \, \fg_0 \, $.  It follows also that the tangent Lie algebra  $ \, \Lie(G_0) = {(\fg_k)}_0 \, $  has the form  $ \, {(\fg_k)}_0 = k \otimes_\Z {(\fg_0)}_{\widetilde{V}} \, $  where  $ \, {(\fg_0)}_{\widetilde{V}} \, $  is the stabilizer of  $ \widetilde{M} $  in  $ \widetilde{V} $:  in turn, this  $ \, {(\fg_0)}_{\widetilde{V}} \, $  depends only on the lattice of weights of  the $\fg_0$-representation
$ \widetilde{V} $  and not on  $ \widetilde{M} $ or
on the choice of a Chevalley basis of  $ \fg_0 $
(see \cite{st} for more details on this classical construction).
 \vskip4pt
   We furthermore require a consistency condition between
$ \fg^\Z $  and  $ G_0 \, $,  as follows.  As the complex Lie algebra
$ \fg $  is simple of classical type, we can fix inside it a
Chevalley basis, as in Sec.~\ref{def_che-sgroups},
call it  $ B $.  Then we assume that
 \vskip2pt
   ---  {\it (a)} \,  $ \; B \cap \fg_0 \, = \, B'_0 \;\; $,
 \vskip2pt
   ---  {\it (b)} \,  $ \,\; \fg^\Z \cap \fg_0 \, = \, {(\fg_0)}_{\widetilde{V}} \;\, $,  $ \;\; \fg^\Z \cap \fg_1 \, = \, \text{\it Span}_\Z\big( B \cap \fg_1 \big) \;\; $.

\vskip11pt

   By  \cite{ccf}, ch.~11,  we have that  $ \, G \subseteq \rGL_{m|n}^k \, $  for suitable $ m $  and  $ n $  and consequently  $ \, \fg_k \subseteq \fgl_{\,m|n}^{\,k} \, $,
%
%
%
 where we denote with  $ \rGL_{m|n}^k $  and  $ \fgl_{\,m|n}^{\,k} \, $  the general linear supergroup and the general linear superalgebra defined over  $ k \, $,  that is  $ \, \rGL_{m|n}^k = \rGL(k^{m|n}) \, $  and  $ \, \fgl_{\,m|n}^{\,k} = \Lie(\rGL_{m|n}^k) \, $,  where  $ k^{m|n} $  is the free  $ k $--supermodule  of dimension  $ m|n $  (see \cite{ccf}, ch.~1, for details).

\smallskip

   Our goal now is to pass from the  $ \bG_0 $--module  $ \, \widetilde{V}_k = k \otimes_\Z \widetilde{V} \, $  to a  $ \bG $--module  $ V_k $  which is obtained as an ``induced representation''  from  $ \bG_0 $  to  $ \bG \, $  (both  $ \widetilde{V}_k $  and  $ V_k $  are  $ k $--modules).  This will be achieved by another ``linearization step'', and an ``induced representation construction'' from $ {\big( \rGL_{m|n}^k \big)}_0 $  to  $ \rGL_{m|n}^k \, $.

\medskip

\begin{remark}
 The results in this section can be easily generalized to the case of  $ k $  a
unital commutative ring, provided we assume  $ G $  to be  {\sl linearizable}.
Notice that this is granted when $k$ is a {\sl field\/}  (see  \cite{ccf},  ch.~11,
and \cite{dg}, ch.~2, for the ordinary setting).  One can check that this is also
granted for  $ k $  a PID and  $ \cO(G) $  a free  $ k $--module.
\end{remark}

We start with a general result on algebraic supergroups, that will be
instrumental to our goal.

\begin{proposition} \label{almost-split}
Let $G$ be an affine algebraic supergroup with  $ \, G \subseteq
\rGL(\cV_k) \, $, for $\cV_k$ a super vector space.
Then we have the following decomposition:
  $$  G \; = \; G_0 \times G_1  \,\; \subseteq \;\,  \rGL(\cV_k)_0 \times \rGL(\cV_k)_1  $$
where  $ G_1 $  is the subscheme defined by  $ \, G_1(A) := G(A) \cap \rGL(\cV_k)_1 \; $.
\end{proposition}

\begin{proof}
Since  $ \, G \subseteq \rGL(\cV_k) \, $,  we have that every  $ \, g \in G(A) \, $
decomposes  in $\rGL(\cV_k)_0 \times \rGL(\cV_k)_1$  uniquely as  $ \, g = g_0 \, g_1 \, $,
with  $ \, g_0 \in \rGL(\cV_k)_0(A) \, $  and  $ \, g_1 \in \rGL(\cV_k)_1(A) \, $ (see \ref{GLdec}).
As  $ \, g_0 = \pi_A \circ g \, $,  where  $ \, \pi_A : A \lra A \big/ J_A \, $,  (as usual $J_A$ denotes the ideal generated by $A_1$ in $A$), we have
that $g_0$ factors via  $ \, \cO(G) \big/ J_{\cO(G)} \, $  and consequently  $ \, g_0 \in
G_0(A) \, $,  \, from which  $ \, g_1 = g_0^{-1} \, g \in G(A) \, $.  Therefore we have
the result.
\end{proof}

\begin{definition}
 With notation as above, let  $ \widetilde{V}_k^{\,*} $  be the  $ G_0 $--module  dual to  $ \widetilde{V}_k \, $.  We define  $ \widetilde{U}_k $  as
  $$  \widetilde{U}_k  \, := \,  \text{\it Ind}_{\,\bG_0}^{\;{(\rGL^k_{m|n})}_0} \big(\, \widetilde{V}_k^{\,*} \big)  $$
i.e.~$ \widetilde{U}_k $  is the  $ \big( \rGL_{m|n}^k \big)_0 \, $--module
%
%
induced from the  $ \bG_0 $--module  $ \widetilde{V}_k^{\,*} \, $.
\end{definition}

\medskip

   Let  $ \widetilde{U}_k^{\,*} $  be the  $ {(\rGL^k_{m|n})}_0 $--module  dual to  $ \widetilde{U}_k \, $;  note that, as  $ \, \text{\it Ind}_{\,\bG_0}^{\;{(\rGL^k_{m|n})}_0} \big(\, \widetilde{V}_k^{\,*} \big) \, $ maps onto  $ \widetilde{V}_k^{\,*} \, $,  we have that  $ \, \widetilde{V}_k \cong \widetilde{V}_k^{\,**} \, $  embeds into  $ \widetilde{U}_k^{\,*} \, $,  i.e.~the latter contains as a  $ \bG_0 $--submodule  an isomorphic copy of  $ \widetilde{V}_k \, $.

\smallskip

   As  $ \widetilde{U}_k^{\,*} $  is a  $ {\big( \rGL^k_{m|n} \big)}_{0} $--module,  it is also a module for the algebra of distributions on  $ {\big( \rGL^k_{m|n} \big)}_{0} \, $,  which identifies with  $ \, \U_k\big(\! {\big( \rgl_{\,m|n} \big)}_0 \big) \! := k \otimes_\Z \U_\Z\big(\! {\big( \rgl_{\,m|n} \big)}_0 \big) $,  the classical Kostant algebra of  $ \, \Lie \big(\! {\big( \rGL^k_{m|n} \big)}_{0} \big) = {\big( \rgl_{\,m|n}^{\,k} \big)}_0 \, $  (cf., for instance,  \cite{ja1}, \S~I.7).  So  $ \widetilde{U}_k^{\,*} $  is a  $ \U_k\big( {\big( \rgl_{\,m|n} \big)}_0 \,\big) $--module,  and we can perform on it the induction from  $ \U_k\big( {\big( \rgl_{\,m|n} \big)}_0 \,\big) $  to  $ \U_k\big( \rgl_{\,m|n} \big) \, $:  this yields next relevant object:

\smallskip

\begin{definition}
  $$  W_k  \; := \;  \text{\it Ind}_{\;\U_k({( \rgl_{\,m|n} )}_0 )}^{\,\;\U_k( \rgl_{\,m|n} )} \big(\, \widetilde{U}_k^{\,*} \big)  \; = \;  \U_k\big(\rgl_{\,m|n}\big) \otimes_{{}_{\scriptstyle
\U_k((\rgl_{\,m|n})_0)}} \widetilde{U}_k^{\,*}
$$
%
%
\end{definition}

\smallskip

\begin{proposition}
Let the notation be as above. Then  $ \, W_k \, $  has a natural structure of  $ \, \rGL^k_{m|n} $--module  and of  $ \, \bG $--module.
\end{proposition}

\begin{proof}
 Clearly, if  $ W_k $  is a  $ \rGL^k_{m|n} $--module  then it is a  $ \bG $--module as well, since  $ \bG $  is a closed subsupergroup of  $ \rGL^k_{m|n} \, $. Let now $\rho$ be the representation map of $ \fgl_{\,m|n}^{\,k} $  into  $ \mathrm{End}(W_k) $  and  $ \sigma $  the representation map of $ (\rGL^k_{m|n})_0 $  into  $ {\mathrm{Aut}}(W_k) \, $.  To give  $ W_k $  a  $ \rGL^k_{m|n} $--module structure, in view of Proposition \ref{GLdec} we need to extend $\sigma$ by specifying the images of all the elements  $ \, I + \theta X \, $  in  $ (\rGL^k_{m|n})_1(A) \, $,  of course in a way compatible with respect to the images of the elements in  $ (\rGL^k_{m|n})_0 \, $.  Let us define
  $$  \sigma( I + \theta X ).w  \, = \, w + \theta \rho(X) w   \eqno \forall \;\; w \in W_k  \qquad  $$
 We leave to the reader the check that this definition is compatible
with the one on  $ (\rGL^k_{m|n})_0 \, $.  This is essentially a consequence
of the fact that  $ \, d\sigma_0 = \rho_0 \, $,  where $\sigma_0$ and $\rho_0$ are
the even parts of the representations $\sigma$ and  $ \rho \, $.
                                            \par
   From another point of view, note that our definition of  $ \, \sigma( I + \theta X ) \, $ is exactly the one giving the unique action of  $ \rGL^k_{m|n} $  on  $ W_k \, $,  induced by restriction of the action of  $ \,\rGL^k_{m|n} \, $,  extending to the action of  $ \rgl_{\,m|n}^{\,k} $  (here we just need to recall that $ \rGL^k_{m|n} $  is naturally embedded into  $ \rgl_{\,m|n}^{\,k} \, $).  In particular, an action of  $ \rGL^k_{m|n} $  on  $ W_k $ with such properties exists, it is unique and it is given exactly by the formula above.
\end{proof}

\smallskip

   Now comes the main result of this subsection.

\medskip

\begin{theorem} \label{embed}
   Let the notation be as above.
 \vskip3pt
   (a) \;  The subspace
  $$  V_k  \; := \;  \U_k(\fg) \otimes_{\U((\fg_k)_0)} \widetilde{V}_k
\, \subseteq \; W_k \, $$
is a rational faithful finite dimensional  $ \, \bG$--module,  and  $ \, G $  embeds into  $ \rGL(V_k) $  as a  {\sl closed}  subsupergroup.
 \vskip3pt
   (b) \;  There exists a Chevalley supergroup  $ \, \bG_V $  such that  $ \, \bG_V \subseteq \rGL(V_k) \, $  and  $ \, \Lie(\bG_V) = \fg_k \; $.  In other words, both  $ \, \bG $  and the Chevalley supergroup  $ \, \bG_V $  embed into the same general linear supergroup  $ \, \rGL(V_k) $  and have the same Lie superalgebra.  Moreover  $ \, \bG_0 = (\bG_V)_0 \; $.
\end{theorem}

\begin{proof}
 First of all, note that by  Remark \ref{super-dist}  we have that  $ \, \U_k(\fg) \subseteq \U_k\big(\fgl_{\,m|n}\big) \, $,  hence  $ V_k $  is a well-defined subspace of  $ W_k \, $:  then by construction, it is also clear that the former is a  $ \bG $--submodule  of the latter.

 Since  $ \widetilde{V}_k $  is rational and faithful as a  $ \bG_0 $--module,  $ V_k $  in turn is rational and faithful as a  $ \bG $--module.  This happens because  $ G $  acts on  $ W_k $  leaving  $ V_k $  invariant.  This is a straightforward application of Proposition \ref{almost-split}.
In particular,  $ \bG $  embeds as a closed subsupergroup inside
$ \, \rGL(V_k) \, $.
                                                      \par
   Now let  $ \widetilde{M} $  be an admissible lattice   --- in the complex  $ \fg_0 $--module  $ \widetilde{V} $  ---   used to construct  $ G_0 $  via a Chevalley construction.  Then we see at
once that  $ \, M := \U_\Z(\fg) \otimes_{{}_{\scriptstyle \U_\Z(\fg_0)}} \! \widetilde{M} \, $
is an admissible lattice for the (rational, faithful) complex  $ \fg $--module  $ \, V := \U_\C(\fg) \otimes_{\U_\C(\fg_0)} \! \widetilde{V} \, $,  which is also finite dimensional because  $ \U_\C(\fg) $  is free of finite rank as a  $ \U_\C(\fg_0) $--module  (cf.~\cite{fg1}, sec.~4).
                                                \par
   Altogether, the above means that we can use  $ V $  and its lattice  $ M $  to construct a Chevalley supergroup  $ G_V $  over  $ k \, $,  realized as a closed subsupergroup of  $ \rGL(V_k) \, $.  As the faithful action of  $ \fg_0 $  onto  $ \widetilde{V} $  yields an embedding of  $ \bG_{V,0} $  into  $ \rGL\big(V_k) \, $,  the restriction to  $ \fg_0 $  of the (faithful) action of  $ \fg $  onto  $ V $  yields an embedding of  $ \bG_{V,0} $  into  $ \rGL\big(V_k) \, $.  By construction   --- including the fact that
 $ \; V_k = \U_k(\fg) \otimes_{\U((\fg_k)_0)} \! \widetilde{V}_k = \bigwedge \! {(\fg_k)}_1 \! \otimes_k \! \widetilde{V}_k \; $
as a  $ \fg_0 $--module  is just  $ \, \widetilde{V}_k^{\,\oplus r} \, $  for  $ \, r := \text{\it rank}_{\,\U((\fg_k)_0)}(\U_k(\fg)) $  ---   the  $ \fg_0 $--action  on  $ V $  is just an  $ r $--fold  diagonalization of the  $ \fg_0 $--action  on  $ \widetilde{V} \, $:  as a consequence, the embedded copy of  $ \bG_{V,0} $  inside  $ V_k $  is just an  $ r $--fold  diagonalized copy of the group obtained from the  $ \fg_0 $--action  on  $ \widetilde{V} $  via the Chevalley construction.  
Hence $ \, \bG_{V,0} = G_0 \, $  inside  $ \rGL(V_k) \, $.
\end{proof}

\medskip

\medskip

\subsection{$G$ as a Chevalley supergroup}

We want to show that $G$ and $G_V$ are isomorphic. Since we shall
make use of the fact that their Lie superalgebras are isomorphic,
we need to make some observations on the differentials.

\begin{lemma}  \label{eqs-diff}
Let  $ \, f \in \cO\big(\rGL(V_k)\big) \, $  and let
$ \, X \in \fgl_1(V_k)(A) \, $,  $ \, A \in \salg \, $
with as usual  $ \, \fgl(V_k) = \Lie(\rGL(V_k)) \, $.  Then
$$
f(1+\theta X) \, = \, f(1) + (df)_1 \theta X  \qquad \qquad \forall
\;\; \theta \in A_1
$$
\end{lemma}

\begin{proof}
Clearly it is enough to check this for a monomial
$ \, f = x_{i_1j_1} \cdots x_{i_rj_r} \, $,  \, where $x_{ij}$ denotes an even or
odd generator of  $ \, \cO\big(\rGL(V_k)\big) \, $.  Notice that
the case of  $ \, f = x_{ij} \, $  is  true:
$x_{ij} (1+\theta X)$ $\; = \;$
$x_{ij}(1) + x_{ij}(\theta X) \; = \; x_{ij}(1) + (dx_{ij})_1 \theta X$.
   The general case reads
$$
 \displaylines{
 (x_{i_1j_1} \dots x_{i_rj_r})(1+\theta X)  \,\; = \;\,
  x_{i_1j_1} (1+\theta X) \cdots x_{i_rj_r}(1+\theta X)  \,\; =   \hfill  \cr
  = \;\, x_{i_1j_1} (1) \cdots x_{i_rj_r}(1) +
x_{i_1j_1} (\theta X) x_{i_2j_2} (1) \cdots  x_{i_rj_r}(1) \, +   \qquad  \cr
  \hfill + \, x_{i_1j_1} (1) x_{i_2j_2} (\theta X) \cdots x_{i_rj_r}(1) + \, x_{i_1j_1} (1) \cdots x_{i_{r-1}j_{r-1}}(1) x_{i_rj_r}(\theta X)  \,\; =   \cr
  \hfill   = \;\, 1 + d(x_{i_1j_1} \cdots x_{i_rj_r})_1(\theta X)  }
$$
which gives what we wanted.
\end{proof}

\medskip

\begin{lemma}
Let the notation be as above. Then  $ \, \bG_V \subseteq \bG \, $,  in other
words  $ \, \bG_V(A) \subseteq \bG(A) \, $  for all  $ A \in \salg \, $.
\end{lemma}

\begin{proof}
 As  $ \bG_V $  is a {\sl closed} subscheme of  $ \rGL(V_k) \, $  (by Theorem
\ref{GVclosed}), an element  $ \, z \in \bG_V(A) \subseteq \rGL(V_k)(A) \, $
corresponds to a morphism  $ \, z : \cO(\rGL(V_k)) \lra A \, $
factoring through  $ I_{\bG_V} \, $,  the ideal defining  $ G_V $
in  $ \cO(\rGL(V_k)) \, $,  that is  $ \, z : \cO(\rGL(V_k))\big/I_{\bG_V}
= \cO(\bG_V) \lra A \, $  (by an abuse of notation we use the same letter).
Hence to prove that  $ \, z \in \bG(A) \, $  we need to show that $z$ factors
also via the ideal $I_\bG$ of  $ \cO(\bG) \, $, which is also closed in $\rGL(V_k)$
(see Theorem \ref{embed}).

\smallskip

If  $ \, z \in (\bG_{V,0})(A) \subseteq \rGL(V_k)_0(A) \, $, then there is nothing to
prove, since  $ \, G_0 = G_{V,0} \, $,  so we assume
$ \, z \in \bG_{V,1}^<(A) \, $ (refer to \ref{mainreschev} for the notation).
It is not restrictive
to assume  $ \, z = 1 + \theta X \, $ for a suitable  $ \, X \in
\fg_1 $ and  $ \, \theta \in A_1 \, $, since such $z$'s together with $G_{V,0}$
generate $G_V(A)$ as an abstract group.
Now let  $ \, f \in I_\bG \, $:  we need to prove that
$$
z(f) \, = \, (1+\theta X)(f) \, = \, f(1+\theta X) \, = \, 0
$$
By the previous lemma we have
$$
f(1+\theta X) \, = \, f(1) + (df)_1 \theta X
$$
Certainly  $ \, f(1) = 0 \, $  because the identity is a topological point
belonging to both $\bG$ and  $ \bG_V \, $.  Moreover,  $ \, (df)_1X = 0 \, $
because of  Proposition 10.6.15 in \cite{ccf},
since $X$ is in the tangent space at the identity
to both supergroups  $ \bG $  and  $ \bG_V \, $.
\end{proof}

\medskip

\begin{lemma}  \label{isovar}
 Let  $ X $  and  $ Y $  two smooth superschemes  (cf.~\cite{fi})  globally split and
such that:
\\
1.  $ \; X \subseteq Y \, $,  $ \; |X| = |Y| \; $;
\\
2.  $ \; T_x X = T_x Y \; $  for all  $ \, x \in |X| \; $.
\\
Then  $ \; X = Y \; $.
\end{lemma}

\begin{proof}
We have a morphism of superschemes
given by the inclusion $X \hookrightarrow Y$.
In order to prove this is an isomorphism it is enough to
verify this on the stalks of the structure sheaves.
The inclusion induces a surjective morphism on the
sheaves, hence we have
$ \; \cO_{Y,y} \twoheadrightarrow \cO_{X,x} \; $.
Since both $X$ and $Y$ are globally split and smooth,
taking completions we have that  $ \, \cO_{X,x} \subseteq \widehat{\cO_{X,x}} \, $  and
$ \, \cO_{Y,y} \subseteq \widehat{\cO_{Y,y}} \, $;  moreover, we can write
the following commutative diagram:
$$
\begin{array}{ccc}
\cO_{Y,x} & \relbar\joinrel\twoheadrightarrow & \cO_{X,x} \\
\downarrow & & \downarrow \\
\widehat{\cO_{Y,x}} & \longrightarrow & \widehat{\cO_{X,x}}
\end{array}
$$
The arrow
  $ \; \widehat{\cO_{Y,x}} \longrightarrow \widehat{\cO_{X,x}} \; $  is
an isomorphism, since both $X$ and $Y$ are smooth and they have
the same tangent space. Hence we have that also the arrow
$ \; \cO_{Y,x} \relbar\joinrel\twoheadrightarrow \cO_{X,x} \; $
is an isomorphism.
\end{proof}

\smallskip

   We are eventually ready for our main result:

\medskip

\begin{theorem} \label{mainresult}
Let  $ \, G $  be an affine algebraic supergroup scheme over the field  $ k \, $,
with  $\bG_0$  being $k$--split,  whose Lie superalgebra $\fg$  is a
$k$--form of a complex Lie superalgebra of classical type.
Then there exists a Chevalley supergroup  $ G_V $  such that  $ \; G_V \cong G \; $.
\end{theorem}

\begin{proof}
Both $G$ and $G_V$ described in the previous propositions
embed into the same $\rGL(V_k)$ and decompose inside the latter
as  $ \, G = G_0 \times G_1 \, $  and  $ \, G_V = G_{V,0} \times G_{V,1} \, $,  with
$ \, G_0 = G_{V,0} \, $.
                                                        \par
   By the previous analysis, we are now left with the following situation:
$ \, G_V \subseteq G \subseteq \rGL(V_k) \, $,  $ \, G_0 = G_{V,0} \, $  and
$ \, T_1{G_V} = T_1{G} \, $.  Actually this happens for all points, not just
the identity, so that  $ \, T_x{G_V} = T_x{G} \, $  for all
$ \, x \in |G| = |G_V |\, $  (notation of ch.~10, sec.~4, in \cite{ccf}).
Then by the lemma \ref{isovar} we have the result, since both
$G$ and $G_V$ are globally split (cf.~\cite{ms1}) and smooth
(since  $ \, G_{V,0} = G_0 \, $  is smooth).
\end{proof}

\smallskip

\begin{observation}
We want to remark that Theorem \ref{mainresult} can be applied in
a different setting, that can be useful for the applications. Assume
$G$ to be a smooth affine algebraic supergroup scheme over a field  $ k \, $:  then
$G$ is a closed subsupergroup scheme in some  $\rGL(V_k)$   --- see \cite{ccf},
ch.~11.  Assume now that $V$ is a suitable representation of
a complex Lie superalgebra  $ \fg \, $,  such that we can construct
the Chevalley supergroup $G_V$ according to the recipe described in
sec.~\ref{chevalley}.  In \cite{fg2} we have shown that such recipe
can be suitably generalized to include Lie superalgebras not of classical
type, for instance the Heisenberg superalgebra.
Assume furtherly that  $ \, G_0 = G_{V,0} \, $  and that
$ \, \Lie(G) = \Lie(G_V) \, $,
in other words $G$ and $G_V$ have the same underlying classical
group scheme and have the same Lie superalgebra.
Then, one can show easily following the arguments
in Theorem \ref{mainresult} that  $ \, G \cong G_V \, $,  that is, our smooth
affine algebraic supergroup $G$ can be realized via the Chevalley
supergroup construction.
\end{observation}

\bigskip

\subsection{\!\!Chevalley Supergroups and Super Harish-Chandra pairs}
\label{shcp-sec}
 In super Lie theory there is an equivalence of categories between
the category of Lie supergroups and the category of Super Harish-Chandra
pairs (SHCP), that is the category consisting of pairs $(G_0, \fg)$, where
$G_0$ is an ordinary real or complex Lie group and
$\fg$ is a real or complex Lie superalgebra
with $Lie(G_0)=\fg_0$ and there is an action of $G_0$ on $\fg$
corresponding to the adjoint action when restricted to $\fg_0$.
Morphisms of SHCP's are defined in a natural way and one can
show a bijective functorial correspondence between the objects and
the morphisms of the given two categories, hence
realizing the equivalence of categories mentioned above (a full
account of the theory is found for example in \cite{ccf},
where the origins of this theory are carefully discussed and references
are given).

\medskip

   A natural question is whether it is possible to extend the theory of SHCP's
to the category of algebraic supergroups.
                                                       \par
   When the algebraic supergroups are over fields of characteristic zero,
the problem has been already treated and solved in  \cite{cf}:  this applies
differential techniques, which cannot be employed instead for arbitrary
characteristic.
                                                       \par
   Instead, more general results are obtained in  \cite{ms2},  using a different
approach, rather closer to the standard one in use for studying algebraic groups
in positive characteristic.  Roughly, one considers a dual version of SHCP where
the first item of the pair is no longer a (classical) algebraic group but a
``hyperalgebra'' instead.  Indeed (still very roughly speaking) if one starts
with an algebraic supergroup  $ G \, $,  then in the corresponding SHCP in the
sense of  \cite{ms2}  the even subgroup  $ G_0 $  is replaced by the (classical)
distribution algebra of  $ G_0 \, $,  the ``correct'' tool for studying  $ G_0 $
in infinitesimal terms.

\medskip

   In the special case of Chevalley supergroups, we can directly prove a certain
equivalence of categories based on the theory developed so far here and in  \cite{fg1}.
As any Chevalley supergroup is built by means of a ``distribution superalgebra'' (namely
the Kostant  $ \Z $--form)  this result is fully consistent with those in  \cite{ms2}.

\medskip

\begin{definition}
 Let  $ k $  be an arbitrary field such that  $ \, \mathrm{char}(k) \neq 2, 3 \, $. We say that $(G_0, \fg)$ is \textit{Chevalley Super Harish-Chandra Pair} (CSHCP), if \\
   \indent   {\it (1)} \;  $G_0$ is an ordinary Chevalley group over  $ k \, $; \\
   \indent   {\it (2)} \;  $\fg$ is a Lie superalgebra of classical type, with  $ \fg_0 = \text{\sl Lie}(G_0) \, $; \\
   \indent   {\it (3)} \;  there is a well defined action, called the \textit{adjoint action} (with a slight abuse of notation) of $G_0$ on $\fg\,$, reducing to the adjoint action on $\fg_0\,$.

\medskip

A \textit{morphism} $(\rho_0, \psi):(G_0,\fg) \lra (H_0,\fh)$ of CSHCPs
consists of a morphism $\rho_0: G_0 \lra H_0$ of algebraic groups
and a morphism $\psi: \fg \lra \fh$ intertwining the adjoint action
of $G_0$ and $H_0$.

\medskip

We shall denote the category of CSHCP with $\CSHCP$.
\end{definition}

\medskip

\begin{proposition}
 There is a unique Chevalley supergroup associated to a given CSHCP.
\end{proposition}

\begin{proof}
Given a CSHCP the recipe given in \cite{fg1} allows us to
produce a Chevalley supergroup associated with it. Section 5.4
in \cite{fg1} proves uniqueness.
\end{proof}

\medskip

We now define $\chesgrp$ to be the category of algebraic supergroups
satisfying the hypothesis carefully detailed at the beginning of
section \ref{unique}.  It is very clear that given $G \in \chesgrp$
there is a unique CSHCP associated with it.  Next theorem
establishes an equivalence of categories.

\medskip

\begin{theorem}
 There exists an equivalence of categories between $\CSHCP$ and $\chesgrp$
%
%
\end{theorem}

\begin{proof}
 The bijective correspondence on the objects is clear, as it is for the morphisms.
\end{proof}

\medskip

\bigskip

\appendix

\section{Chevalley basis}

In this appendix we quickly recall the definition of Chevalley basis
(see \cite{fg1} for more details).

\medskip

Assume $\fg$ to be a Lie superalgebra of classical type different
from $A(1,1)$, $P(3)$, $Q(n)$ and $D(2,1;a)$, $a \notin \Z$.
We prefer to leave out these cases to simplify our definitions,
for a complete treatment see \cite{fg1}.

\medskip

   Let us fix a Cartan subalgebra $\fh$ of $\fg \, $:  its adjoint action
gives the {\sl root space} decomposition of $\fg$
$$
\fg \; = \; \fh \, \oplus \, {\textstyle \bigoplus_{\alpha \in \Delta}} \fg_\alpha
$$
where  $ \, \Delta = \Delta_0 \cup \Delta_1 \, $  is the root system, with
$$
\begin{array}{c}
\Delta_0  \, := \,  \big\{\, \alpha \in \fh^* \setminus \{0\} \;\big|\;
\fg_\alpha \cap \fg_0 \not= \{0\} \big\}  \, = \,
\{\,\hbox{{\sl even roots\/}  of  $ \fg $} \,\}.  \\ \\
\Delta_1  \, := \,  \big\{\, \alpha \in \fh^* \;\big|\;
\fg_\alpha \cap \fg_1 \not= \{0\} \big\}  \, = \,
\{\,\hbox{{\sl odd roots\/}  of  $ \fg $} \,\}.
\end{array}
$$
If we fix a simple root system (see \cite{ka} for its definition)
the root system splits into positive and negative roots, exactly as
in the ordinary setting:
$$
\Delta \, = \, \Delta^+ \,{\textstyle \coprod}\, \Delta^- \;\; ,  \qquad
\Delta_0 \, = \, \Delta_0^+ \,{\textstyle \coprod}\, \Delta_0^- \;\; ,  \qquad
\Delta_1 \, = \, \Delta_1^+ \,{\textstyle \coprod}\, \Delta_1^- \;\; .
$$

\smallskip

   If  $ \fg $  is neither of type  $P(n)$  nor  $Q(n) \, $,
there is an even non-degenerate, invariant bilinear form on
$\fg \, $,  whose restriction to  $ \fh $  is in turn an invariant
bilinear form on  $ \fh \, $.  On the other hand, if $ \fg $
is of type  $ P(n) $  or  $ Q(n) \, $,  then such a form on
$ \fh $  exists because  $ \fg_0 $
is simple (of type  $ A_n $), though it does not come by restricting
an invariant form on the whole $\fg \, $.

\medskip

If $ \, \big( x, y \big) \, $ denotes such form,  we can identify  $ \fh^* $
with  $ \fh $,  via  $ \, H'_\alpha \mapsto \big(H'_\alpha,\ \big) \, $.
We can then transfer $\big( ,  \big)$ to $\fh^*$ in the natural way:
$\big( \alpha, \beta \big) = \big( H'_{\alpha}, H'_{\beta} \big) \; $.
Define  $ \, H_\alpha:= 2{H'_\alpha \over \big(H'_\alpha, H'_\alpha\big)} \, $
when the denominator is non zero.
When $\big(H'_\alpha, H'_\alpha\big)=0$ such renormalization can
be found in detail in \cite{ik}.
We call $H_\alpha$ the \textit{coroot} associated with $\alpha\,$.

\medskip

 \begin{definition}  \label{def_che-bas}
We define a  {\it Chevalley basis\/}  of  a Lie superalgebra $ \fg $
as above any homogeneous  basis
$$
  B  \, = \,  \big\{ H_1 \dots H_l, \, \, X_\alpha, \, \, \alpha \in \Delta \big\}
$$
of $\fg$ as complex vector space, with the following requirements:
 \vskip8pt
 \textit{(a)}
$ \, \big\{ H_1 , \dots , H_\ell \big\} \, $  is a
basis of the complex vector space $\fh\,$.  Moreover
$$
\fh_\Z  \, := \,  \text{\it Span}_{\,\Z} \big\{ H_1 , \dots , H_\ell \big\}
\, = \,  \text{\it Span}_{\,\Z} \big\{ H_\alpha \,\big|\, \alpha \!
\in \! \Delta \}
$$
 \vskip8pt
%
 \textit{(b)}  \hskip4pt   $ \big[ H_i \, , H_j \big] = 0 \, ,   \hskip9pt
 \big[ H_i \, , X_\alpha \big] = \, \alpha(H_i) \, X_\alpha \, ,   \hskip15pt  \forall \; i, j \! \in \! \{ 1, \dots, \ell \,\} \, ,  \; \alpha \! \in \! \Delta \; $;
 \vskip11pt
   \textit{(c)}  \hskip7pt   $ \big[ X_\alpha \, , \, X_{-\alpha} \big]  \, = \,  \sigma_\alpha \, H_\alpha  \hskip25pt  \forall \;\; \alpha \in \Delta \cap (-\Delta) $
 \vskip4pt
\noindent
 with  $ H_\alpha $  suitably defined exactly as in the ordinary setting,
and  $ \; \sigma_\alpha := -1 \; $  if  $ \, \alpha \in \Delta_1^- \, $,  $ \; \sigma_\alpha := 1 \; $  otherwise;
 \vskip13pt
   \textit{(d)}  \quad  $ \, \big[ X_\alpha \, , \, X_\beta \big]  \, = \, c_{\alpha,\beta} \; X_{\alpha + \beta}  \hskip17pt   \forall \;\, \alpha , \beta \in \Delta \, : \, \alpha \not= -\beta \, $,  \, with  $ \, c_{\al,\be} \in \Z \, $.  More precisely,

\begin{itemize}
\item  If  $ \, (\al,\al) \neq 0 \, $,  or  $ \, (\beta,\beta) \neq 0 \, $,  then
$ \, c_{\al,\be} =\pm(r+1) \, $  or (only if  $ \, \fg = P(n) \, $)  $ \, c_{\al,\be} = \pm(r+2) \, $,
where  $ r $  is the length of the  $ \alpha $--string  through  $ \beta \, $.

\item If  $ \, (\al,\al) = 0 = (\beta,\beta) \, $,  then  $ \, c_{\al,\be} = \be(H_\alpha) \, $.
\end{itemize}
\end{definition}

   Notice that this definition clearly extends to direct sums of finitely
many  of the $ \fg $'s under the above hypotheses.

\smallskip

\begin{definition}
 If  $ B $  is a Chevalley basis of  a Lie superalgebra $ \fg \, $ as above, we set
$$
\fg_\Z := \text{span}_\Z\{B\}  \quad  \big(\! \subseteq \fg \,\big)
$$
and we call it the  {\it Chevalley superalgebra\/}  of  $ \fg $.
\end{definition}

\smallskip

   Observe that $ \fg_\Z $  is a Lie superalgebra over  $ \Z \, $ inside  $ \fg \, $.

%
%

%
%
%

\vskip13pt

  $$  \begin{matrix}
   \text{Rita Fioresi} \quad  &  \quad  \text{Fabio Gavarini}  \\
   \text{Dipartimento di Matematica}  \qquad  &  \qquad  \text{Dipartimento di Matematica}  \\
   \text{Universit\`a di Bologna}  \qquad  &  \qquad  \text{Universit\`a di Roma ``Tor Vergata''}  \\
   \text{Piazza di Porta San Donato, 5}  \qquad  &  \qquad  \text{via della ricerca scientifica 1}  \\
   \text{I--40127 Bologna, Italy}  \qquad  &  \qquad  \text{I--0133 Roma, ITALY}  \\
   \text{e-mail:  {\tt fioresi@dm.unibo.it}}  \qquad  &  \qquad  \text{e-mail:  {\tt gavarini@mat.uniroma2.it}}
      \end{matrix}  $$

\end{document}